# Note On Prime Gaps And Very Short Intervals
N. A. Carella, September, 2010

**Abstract:** Assuming the Riemann hypothesis, this article discusses a new elementary argument that seems to prove that the maximal prime gap for a finite sequence of prime numbers $p_1, p_2, \ldots, p_n \leq x$, satisfies $\max\{ d_n = p_{n+1} - p_n : p_n \leq x \} \leq c_0(\log x)^2/\log\log x$. Equivalently, it shows that the very short intervals $(x, x + y]$ contain prime numbers for all real numbers $y \geq c_1(\log x)^2/\log\log x$, where $c_0 > 0$, $c_1 > 0$ are constants, and $x > 0$ is a sufficiently large number.

## 1. INTRODUCTION

For $n \geq 1$, let $p_n$ denotes the $n$th prime in the set of prime numbers $\mathbb{P} = \{ 2, 3, 5, 7, \ldots, \}$ and let $\pi(x) = \#\{ p \leq x : p \text{ prime} \}$ denotes the prime counting function. The *prime gap* $d_n = p_{n+1} - p_n$ is the difference of two consecutive prime numbers, and the maximum prime gap is defined by $d_{\max} = \max_{p_n \leq x}\{ d_n = p_{n+1} - p_n \}$. The average prime gap of the short sequence $p_1, p_2, \ldots, p_n \leq x$ is given by the asymptotic formula $\overline{d}_n = x/\pi(x) = \log x + O(1)$ for a sufficiently large number $x > 0$. The prime gaps and the existence of prime numbers in short intervals, the local distribution, and related concepts, are intensively investigated local properties of the prime numbers, see [RN], [SD] for surveys, and [IK], [GPY], [GV] et cetera for advanced materials.

A resolution of the limit infima of the best parameters $\alpha \geq 0$ and $\beta \geq 0$, for the prime gap inequality

$$p_{n+1} - p_n < p_n^\alpha \log^\beta p_n \tag{1}$$

has been an open problem of much interest for quite sometime. According to the literature, the result below seems to be the sharpest possible that can be derived from the Riemann hypothesis.

***Theorem* 1.** Assume the Riemann hypothesis, then $p_{n+1} - p_n \ll p_n^{1/2} \log p_n$.

In the other direction, the current attempts to show the existence of an infinite sequence of prime numbers of very large gaps is not even close to the conjectured maximal $p_{n+1} - p_n \leq c_0(\log x)^2$. The best so far is the infinite sequence of Westzynthius primes $\{ p_n \in \mathbb{P} : n \in W \}$ of very large gaps

$$p_{n+1} - p_n \geq c_1 \log p_n \log\log p_n \log\log\log\log p_n (\log\log\log p_n)^{-2}, \tag{2}$$

where $c_1 > 0$ is a constant, and $W \subset \mathbb{N}$ is an infinite subset of integers, see Theorems 5 and 6.

A conjecture of Cramer claims that the maximal prime gap satisfies $p_{n+1} - p_n \leq (\log p_n)^2$, see [CR], [GV], and a conjecture of Erdos claims that there exists an infinite sequence of prime numbers of very large gaps $p_{n+1} - p_n \geq c_1(p_n) \log p_n \log\log p_n \log\log\log\log p_n (\log\log\log p_n)^{-2}$, where $c_1(x) \to \infty$ as $x \to \infty$, see [ER], [MP]. It is instructive to compare this result (and the twin prime conjecture) to the conjectured Poisson distribution of the primes gaps, that is,

$$\#\{p_n \leq x : p_{n+m} - p_n > \lambda\} \sim x \sum_{0 \leq k \leq m} \lambda^k e^{-\lambda}/k!,$$

confer, [GH]. Assuming the validity of the Riemann hypothesis, both of these conjectures are very close to the true maximal prime gap $p_{n+1} - p_n = O(\log^2 p_n / \log\log p_n)$. The conjectured interval

$$0.122985 d^{1/2} e^{d^{1/2}} < p_n < 2.096 d e^{d^{1/2}} \tag{3}$$

that contains a pair of prime numbers $q$ and $p$ and the first appearance of the specified prime gap $d = q - p \geq 2$ also seems to be quite close to the true interval, see [NT] for additional information. This is related to a conjectured magnitude of the first appearance of a prime $p \sim \sqrt{d} e^{.5(\log^2 d + 4d)^{1/2}}$ such that the prime gap $d = q - p \geq 2$ with $p < q$. This improves the earlier estimate $p = e^{(1+o(1))d^{1/2}}$ of Shanks in [SK].

The new result, contributed to the literature, is the following.

***Theorem 2.*** (i) If the Riemann hypothesis is true, then $p_{n+1} - p_n = O(\log^2 p_n / \log\log p_n)$.
(ii) The $n$th prime gap satisfies the inequality $p_{n+1} - p_n < O(p_n^{.1559458+\varepsilon})$, $\varepsilon > 0$, unconditionally.

Statement (i) is similar to a result in [SG]. There, assuming the Riemann hypothesis, Selberg proved that $p_{n+1} - p_n \leq f(n)\log^2 p_n$ for almost all $n \geq 1$, where $f(n) \to \infty$ as $n \to \infty$. And statement (ii) is similar to a result in [HN]. There, Harman proved that $p_{n+1} - p_n \leq p_n^{1/10+\varepsilon}$ for almost all $n \geq 1$, where $\varepsilon > 0$ is a small real number.

The implied constant in estimate (i) is less than $2\pi$, see Theorem 11. The numerical data compiled by several authors over the last decades is 100 percent consistent with this estimate. Exempli gratia, as of 2010, the largest maximal prime gap known $1476 = p_{n+1} - p_n \leq 2\pi(\log p_n)^2/\log\log p_n = 2941.0...$ for $p_n = 1425172824437699411$ satisfies it, see the table in Section 4 for other examples.

***Corollary 3.*** Assuming the Riemann hypothesis, for all sufficiently large numbers $x > 0$, and $y \geq c \log^2 x / \log\log x$, some $c > 2\pi$ constant, the very short intervals $(x, x + y]$ contains primes.

**Note:** The *prime existence problem* and the *prime density problem* for very short intervals $(x, x + y]$, $0 < y \leq x^\alpha$, $\alpha > 0$, should not be confused. The latter seems to be a much more difficult problem than the former. The prime existence problem asks for the smallest real number $y > 0$ such that the very short interval $(x, x + y]$ contain primes, id est,

$$\pi(x + y) - \pi(x) > 0.$$

The prime existence problem is equivalent to the maximal prime gap problem. Currently, it is known that $\pi(x + y) - \pi(x) > 0$ for $y > x^{.525}$ unconditionally, or equivalently $p_{n+1} - p_n \leq cp_n^{.525}$, $c > 0$ constant, see [BK]. Furthermore, $\pi(x + y) - \pi(x) > 0$ for almost every $y > x^{1/10+\varepsilon}$ unconditionally, see [HN].

In contrast, the prime density problem asks for the smallest real number $y > 0$ such that the density of prime numbers in the very short intervals $(x, x + y]$ is given by the asymptotic formula

$$\pi(x+y) - \pi(x) \sim y/\log x,$$

see [SG] for a discussion. It is known that $\pi(x+y) - \pi(x) \sim y/\log x$ for $y > x^{7/12}$ unconditionally, see [HB]. But for very short intervals, it is known that $\pi(x+y) - \pi(x) \neq y/\log x$ for $y = \log^B x$, $B > 0$ constant, for infinitely many $x \geq x_0$, see [MH]. Furthermore, $\pi(x+y) - \pi(x) > y/8\log x$ for almost every $y > x^{1/10+\varepsilon}$ unconditionally, see [HN].

Most of the literature on prime gaps, idem quod primes in short intervals, appears under the subject of primes in short intervals. The articles [MW], [GPY], [MJ] and [PT] present surveys and new results on the state of knowledge in the theory of prime numbers and related topics. Other references for other specialized subjects are also given throughout the paper.

The proof of Theorem 2 appears in Section 3. A few of the results required in the proofs of the new results are given in Section 2. Throughout the paper the values of the constants $c$, $c_0$, $c_1$, $c_2$, ..., are local constants, and can vary from result to result.

## 2. FUNDAMENTAL BACKGROUND MATERIALS

The basic principle employed to achieve the new result on the maximal prime gaps and the existence of primes in very short intervals is far simpler than the established theories such as *zero density methods*, see [IK], [IV], [HL], and *sieve methods*, see [BK]. The basic principle is derived from a few elementary results on the theory of prime numbers, and the zeros of the zeta function. Information on these elementary concepts are provided in this Section.

**2.1 Formula for the $n$th Prime.** A formula for the $n$th prime in term of its position in the sequence of prime numbers 2, 3, 5, 7, … is stated here.

*Corollary* **4.** Let $n$ be a large number and let $p_n$ be the $n$th prime. Then the following statements are satisfied.
(i) The $n$th prime satisfies $an \log n \leq p_n \leq bn \log n$ for some constants $a, b > 0$.
(ii) $p_n \sim n \log n$ as the integer $n \to \infty$.

This result is readily deduced from the Prime Number Theorem, confer [EL], [HW, p. 12] and similar sources for related proofs.

The asymptotic formula $p_n = n \log n + o(n \log n)$ probably cannot be used to estimate the $n$th prime gap $d_n = p_{n+1} - p_n$. This obstruction seems to arise from the fact that in the difference of two consecutive primes $p_{n+1} - p_n = \log n + o(n \log n)$, the error term is larger than the main term.

In an arithmetic progression $q_n \equiv a \bmod q$, $\gcd(a, q) = 1$, the $n$th prime $q_n$ in the progression has the asymptotic formula $q_n \sim \varphi(q) n \log n$, see [CL] for a proof. In this case, assuming the generalized Riemann hypothesis, the conjectured prime gaps in arithmetic progression is $q_{n+1} - q_n = O(\varphi(q) \log^2 q_n / \log \log q_n)$.

**2.2 Large Prime Gaps.** The prime gaps can be arbitrarily large as illustrated by the sequence of integers

$$m!, \; m!+1, \; m!+2, \; \ldots, \; m!+m, \; m!+m+1, \; \ldots \; .$$

Assuming $p_n = m! + 1$ is prime, the next potential prime is $p_{n+1} = m! + m + 1$. Ergo, the minimal prime gap in this sequence of numbers satisfies $p_{n+1} - p_n = m \geq \log p_n / \log\log p_n$. Refer to [CT] for numerical data on the primorial primes $p_n = m! \pm 1$.

The most important techniques for constructing an infinite sequence of consecutive prime numbers with very large prime gaps above the average value, viz, $d_n > \bar{d}_n \sim \log p_n$ arise in the proof of the following result, see [MW] for a survey.

***Theorem 5.*** (Westzynthius) Let $x > 1$ be a large real number, and let $p_n \leq x$ be the $n$th prime. Then $p_{n+1} - p_n \geq c \log x \log_3 x (\log_4 x)^{-1}$ for infinitely many $n \geq 1$ as $x \to \infty$. In particular, the limit $\limsup_{n \to \infty} (p_{n+1} - p_n) / \log p_n = \infty$ holds.

There are several ways of proving this result. These techniques use both elementary methods and advanced complicated methods, see [MV, p. 221], and [MP]. The Westzynthius method has been developed by several authors. These authors have achieved concrete formulas and resolved the implied constant, see [ER], [RA], [MP], and [PT].

***Theorem 6.*** ([PT]) There is an infinite sequence of pairs of consecutive primes such that

$$p_{n+1} - p_n \geq c \frac{\log x \log\log x \log\log\log\log x}{(\log\log\log x)^2} \tag{4}$$

for $p_n \leq x$, and $c = 2(e^\gamma + o(1))$ is a constant.

### 2.3 Zeta Zeros Results

The relevant formula for the counting function $N(T)$ for the number of zeros of the zeta $\zeta(s)$ function on the rectangle $R(T) = \{ s \in \mathbb{C} : 0 \leq \Re(s) < 1 \text{ and } 0 \leq Im(s) \leq T \}$ on the critical strip $\{ s \in \mathbb{C} : 0 \leq \Re(s) < 1 \}$ was proposed by Riemann and proved by vonMangoldt.

***Theorem 7.*** (vonMangoldt) Let $N(T) = \#\{ s \in R(T) : \zeta(s) = 0 \}$ be the counting function for the number of complex zeros of height $T \geq 0$. Then

$$N(T) = \frac{T}{2\pi} \log \frac{T}{2\pi} - \frac{T}{2\pi} + O(\log T) + S(T), \tag{5}$$

where $S(t) = \pi^{-1} \arg(\zeta(1/2 + it))$, and $T \geq 1$.

For a detailed proof, see [IV, p. 19], [ES, p. 127], [EL, p. 160] and similar literature.

***Corollary 8.*** The critical zeros $\rho_n = 1/2 + i\gamma_n$ of $\zeta(s)$ satisfy the following.
(i) The imaginary part satisfies $an / \log n \leq \gamma_n \leq bn / \log n$ for some constants $a, b > 0$.
(ii) The imaginary part $\gamma_n \sim 2\pi n / \log n$ as the integer $n \to \infty$.

More details on this result appear in [EL, p. 160], [IV, p. 20], et cetera. This formula is probably not an effective way of estimating the difference $\gamma_{n+k} - \gamma_n = k / \log(n+k) + o(n / \log n)$, since the error term here seems to be larger than the main term.

The study of the distribution of the zero spacings, known as the GUE hypothesis, is considered in [OD], [GN, p. 27] and other.

## 2.4 Prime-Zero Duality Principle
The idea of changing domain from the primes domain to the zeta zeros domain to study the properties of the prime numbers dates back to Riemann as sketched in the *explicit formula*. The *prime-zero duality principle* is a simple transformation from the primes domain to the zeros domain. Information on the zeros of the zeta function (more generally *L*-functions) has dual corresponding information on the prime numbers and vice versa. Accordingly, the duality principle can be used to determine unknown information on the primes from the corresponding known information on the zeros and vice versa.

**Lemma 9.** (Prime-Zero-Duality) For a sufficiently large integer $n \geq 1$, let $p_n$ be the $n$th prime and let $\rho_n = \beta_n + i\gamma_n$ be the $n$th zero of the zeta function. Then the following asymptotic formulae

(i) $p_n = \dfrac{\gamma_n \log^2 n}{2\pi} + o(\gamma_n \log^2 n)$,  (ii) $\gamma_n = \dfrac{2\pi p_n}{\log^2 n} + o\left(\dfrac{p}{\log^2 n}\right)$,

hold.

*Proof*: By Corollaries 4 and 8, the asymptotic ratio

$$\frac{p_n}{\gamma_n} = \frac{n \log n + o(n \log n)}{c_0 n / \log n + o(n / \log n)} = c_1 \log^2 n + o(\log^2 n), \qquad (6)$$

holds for large integers $n \geq 1$, where $c_0 = 2\pi$, $c_1 = 1/2\pi$. This leads to the asymptotic expression $p_n = (c_1 \log^2 n + o(\log^2 n))\gamma_n$, which proves (i). The proof of statement (ii) uses similar analysis. ∎

As shown before, a direct application of formula 8-ii is probably not an effective way of estimating the zero spacing, but there are other means of accomplishing this task.

**Theorem 10.** ([IV, p. 261]) For any $\varepsilon > 0$ and $n \geq n_0(\varepsilon)$, the $n$th critical zero spacing satisfies the inequality

$$\gamma_{n+1} - \gamma_n < \gamma_n^{\theta + \varepsilon}$$

where $\theta = .1559458...$ .

This is an unconditional result for the real parts $\gamma = \Re(\rho)$ of the critical zeros $\rho = \sigma + i\gamma$ of the zeta function. A significantly sharper but conditional result is also known.

**Theorem 11.** ([CG]) Assume the Riemann hypothesis. Then $\gamma_{n+1} - \gamma_n \leq \dfrac{\pi}{\log \log \gamma_n}(1 + o(1))$.

## 3. MAIN RESULT
The first subsection is a short introduction to the theory of primes in short intervals. This is followed by a proof of the main result.

### 3.1 Primes in Very Short Intervals
The earliest result on primes in short intervals was achieved by Tchebychev, it deals with intervals of the form $(x, x + y]$ with $y = cx$, $c \approx 2$. The analysis was further developed to $y \asymp x$, see [NZ, p. 224], [HW, p. 455]. The special case of $c = 2$ is widely known as Bertrand's Postulate. The Gauss form $\pi(x) = x/\log x + o(x/\log x)$ of the Prime Number Theorem implies the existence of primes in the short intervals $(x, x + y]$ with $y \asymp x/\log^\varepsilon x$, $\varepsilon < 1$. And the DelaValle Poussin form $\pi(x) = li(x) + O(xe^{-c\sqrt{\log x}})$ implies the existence of primes in the short intervals $(x, x + y]$ of subexponential sizes $y = xe^{-c(\log x)^{1/2-\varepsilon}}$ with $c > 0$ and $\varepsilon > 0$ constants.

The existence of primes in the short intervals $(x, x + y]$ of exponential size $y = x^{1-\varepsilon}$, $1/3300 < \varepsilon < 1$, or better, was achieved by Hoheisel. The analysis, often called the *zero density method*, revolves around the vonMangoldt explicit formula and related concepts, see [IK], [IV], [MV] and similar references. The current record is approximately $y = x^{1/2+\varepsilon}$, $1/12 < \varepsilon < 1$, see below.

***Theorem 12.*** ([HL]) Let $\theta > 7/12$. Then the following holds.
(i) The interval $(x, x + y]$ contains $y(1 + o(1))/\log x$ primes for $y \geq x^\theta$.
(ii) Let $p_n$ be the $n$th prime, then $p_{n+1} - p_n \leq p_n^\theta$.

The sieve methods are analytical methods based on the modern theory of the sieve of Eratosthenes. The most recent result on the application of sieve methods to the theory of primes in short intervals is the following.

***Theorem 13.*** ([BK]) For all large $x$, the interval $(x - x^{1/2+1/40}, x]$ contains primes numbers.

The Riemann hypothesis limits the size of the short intervals that can be analyzed using the explicit formula and the zero density method. This limitation is probably true for the sieve methods too. The best possible result under this hypothesis is stated here.

***Theorem 14.*** (vonKoch) If the Riemann hypothesis holds, then
(i) There are at least $cx^{1/2} \log x$ primes in the interval $[x, x + x^{1/2+\varepsilon}]$, $c > 0$ constant.
(ii) The prime gap is of order $p_{n+1} - p_n = O(p_n^{1/2} \log^2 p_n)$.
(iii) $\psi(x+y) - \psi(x) = y + O(x^{1/2} \log^2 x)$.

Proof: These are derived from the integral $\pi(x+y) - \pi(x) = \int_x^{x+y} (\log t)^{-1} dt + O(x^{1/2} \log^2 x)$, where $y = x^{1/2+\varepsilon}$, see [NW, p. 245]. For a different proof of (ii), see [IV, p. 321]. ∎

**3.2 Main Result**
The average prime gap of the short sequence $p_1, p_2, \ldots, p_n$ of prime numbers is given by the asymptotic formula $\bar{d}_n = p_n/\pi(p_n) = \log p_n + O(1)$, and the average zero spacing of the short sequence $\gamma_1, \gamma_2, \ldots, \gamma_n$ is given by the asymptotic formula $\bar{\delta}_n = \gamma_n/N(\gamma_n) = 1/\log \gamma_n + O(\log \gamma_n / \gamma_n)$. These statistics are inversely proportional. And the Prime-Zero Duality formula, (Lemma 9), seems to imply that there is a one-to-one correspondence between the sequence of prime numbers $\{p_n : n \geq 1\}$ and the sequence of real parts $\{\gamma_n : n \geq 1\}$ of the zeta zeros.

***Theorem 2.*** (i) If the Riemann hypothesis is true, then $p_{n+1} - p_n = O(\log^2 p_n / \log\log p_n)$.
(ii) The $n$th prime gap satisfies the inequality $p_{n+1} - p_n < O(p_n^{.1559458+\varepsilon})$, $\varepsilon > 0$, unconditionally.

Proof of (i): Fix a small integer $k < n$. Applying Corollary 8 to compute the difference of two distinct zeros returns

$$\gamma_{n+k} - \gamma_n = \left( \frac{2\pi \, p_{n+k}}{\log^2(n+k)} + o\left(\frac{p_{n+k}}{\log^2(n+k)}\right) \right) - \left( \frac{2\pi \, p_n}{\log^2 n} + o\left(\frac{p_n}{\log^2 n}\right) \right)$$
$$= (p_{n+k} - p_n)\frac{2\pi}{\log^2 n} + o\left(\frac{p_{n+k}}{\log^2 n}\right) - o\left(\frac{p_n}{\log^2 n}\right), \tag{7}$$

where the implied constants are distinct, and the effect of the constant $k \geq 1$ is absorbed in the error term. Now assuming the Riemann hypothesis, the $n$th zero spacing satisfies the inequalities $0 < \gamma_{n+k} - \gamma_n \leq c_0 / \log\log \gamma_n$, confer Theorem 11. Since $n = c \log \gamma_n$, $c > 0$ constant, this is rewritten as a function of $n \geq 1$ as

$$0 < \gamma_{n+k} - \gamma_n \leq \frac{c_1}{\log\log n}, \tag{8}$$

where $c_0 > 0$, $c_1 > 0$, $c_2 > 0$, …, are constants. Combining (7) and (8) returns

$$0 < (p_{n+k} - p_n)\frac{2\pi}{\log^2 n} + o\left(\frac{p_{n+k}}{\log^2 n}\right) - o\left(\frac{p_n}{\log^2 n}\right) = \gamma_{n+k} - \gamma_n \leq \frac{c_1}{\log\log n}. \tag{9}$$

Simplifying (9) yields the equivalent expression

$$0 < p_{n+k} - p_n + (o(p_{n+k}) - o(p_n)) < c_2 \log^2 n / \log\log n. \tag{10}$$

The verification of the claim will be broken up into four separate cases.

**Case 1.** Assume that the error term satisfies $|o(p_{n+1}) - o(p_n)| < c_2 \log^2 n / \log\log n$ for all integers $n \geq n_0$. By the hypothesis on the size of the error term, the main inequality

$$0 < p_{n+1} - p_n + (o(p_{n+1}) - o(p_n)) < c_2 \log^2 n / \log\log n \tag{11}$$

immediately leads to the inequalities $2 \leq p_{n+1} - p_n < 2c_2 \log^2 n / \log\log n$.

**Case 2.** Assume that the error term satisfies $o(p_{n+1}) - o(p_n) > c_2 \log^2 n / \log\log n > 0$ for all integers $n \geq n_0$. Put $o(p_{n+1}) - o(p_n) \geq c_2 \log^{2+\alpha} n / \log\log n > 0$, where $\alpha > 0$ is an arbitrarily small real number. Then the main inequality

$$0 < p_{n+1} - p_n + (o(p_{n+1}) - o(p_n)) < c_2 \log^2 n / \log\log n$$

becomes

$$0 < p_{n+1} - p_n + c_2 \log^{2+\alpha} n / \log\log n \leq p_{n+1} - p_n + (o(p_{n+1}) - o(p_n))$$
$$< c_2 \log^2 n / \log\log n. \tag{12}$$

Clearly, this is a contradiction for all integers $n \geq n_0$.

**Case 3.** Assume that the error term satisfies $o(p_{n+1}) - o(p_n) < -c_2 \log^2 n / \log\log n < 0$ and $p_{n+1} - p_n < c_2 \log^2 n / \log\log n$ for all integers $n \geq n_0$.

Put $o(p_{n+1}) - o(p_n) \leq -c_2 \log^{2+\alpha} n / \log\log n < 0$, where $\alpha > 0$ is an arbitrarily small real number. Then the main inequality

$$0 < p_{n+1} - p_n + (o(p_{n+1}) - o(p_n)) < c_2 \log^2 n / \log\log n \tag{13}$$

becomes

$$\begin{aligned} 0 < p_{n+1} - p_n + (o(p_{n+1}) - o(p_n)) &\leq p_{n+1} - p_n - c_2 \log^{2+\alpha} n / \log\log n \\ &< c_2 \log^2 n / \log\log n. \end{aligned} \tag{14}$$

Clearly, this is a contradiction for all $2 \leq p_{n+1} - p_n < c_2 \log^2 n / \log\log n$ and all integers $n \geq n_0$.

**Case 4.** Assume that the error term satisfies $o(p_{n+1}) - o(p_n) < -c_2 \log^2 n / \log\log n < 0$ and $p_{n+1} - p_n > c_2 \log^2 n / \log\log n > 0$ for all integers $n \geq n_0$.

Put $o(p_{n+1}) - o(p_n) \leq -c_2 \log^{2+\alpha} n / \log\log n < 0$, and $p_{n+1} - p_n \geq c_2 \log^{2+\beta} n / \log\log n > 0$, where $\alpha > 0$ and $\beta > 0$ are arbitrarily small real numbers. Then the main inequality is rewritten as

$$\begin{aligned} 0 < c_2 \log^{2+\beta} n / \log\log n + (o(p_{n+1}) - o(p_n)) &\leq p_{n+1} - p_n + (o(p_{n+1}) - o(p_n)) \\ &< c_2 \log^2 n / \log\log n. \end{aligned} \tag{15}$$

Further simplification yields

$$0 < 1 + \frac{o(p_{n+1}) - o(p_n)}{c_2 \log^{2+\beta} n / \log\log n} < \frac{1}{\log^\beta n}. \tag{16}$$

If $\alpha > \beta > 0$, and $o(p_{n+1}) - o(p_n) \leq -c_2 \log^{2+\alpha} n / \log\log n < 0$ for all integers $n \geq n_0$, then

$$0 < 1 + \frac{o(p_{n+1}) - o(p_n)}{c_2 \log^{2+\beta} n / \log\log n} < \frac{-1}{2} < \frac{1}{\log^\beta n} \tag{17}$$

for all integers $n \geq n_0$. And if $0 < \alpha < \beta$, and $o(p_{n+1}) - o(p_n) \leq -c_2 \log^{2+\alpha} n / \log\log n < 0$ for all integers $n \geq n_0$, then

$$0 < \frac{1}{2} < 1 + \frac{o(p_{n+1}) - o(p_n)}{c_2 \log^{2+\beta} n / \log\log n} < \frac{1}{\log^\beta n} \tag{18}$$

for all integers $n \geq n_0$. In either situation (15) or (16), this is a contradiction for all pair $\alpha > 0$, $\beta > 0$ such that $\alpha \neq \beta$, and large $n$.

The value of the parameters $\alpha = \beta > 0$ does not occur since $\gamma_{n+1} - \gamma_n > 0$. Therefore, Case 1 is the only possibility. The verification of statement (ii) is similar as above but uses Theorem 10 instead. ∎

## 4. NUMERICAL DATA

Significant amount of data have been complied by several authors over the past decades, see [NT], [NL], and [TS] for full details. A summary of some of these data is provided in a table for comparative study. The fourth column records the estimated maximal prime gap $p_{n+1} - p_n = 2\pi \log^2 p_n / \log\log p_n = U(n)$.

As shown before, the largest maximal prime gap known $1476 = p_{n+1} - p_n \leq 2\pi(\log p_n)^2/\log\log p_n = 2941.0...$ for $p_n = 1425172824437699411$, see the table. It is very interesting to note that the interlacing composites

$p_n = 1425172824437699411$,
$p_{n+1} = 1425172824437699412 = 2^2 \cdot 3^3 \cdot 43 \cdot 601 \cdot 510623560373$,
$p_{n+2} = 1425172824437699413 = 17 \cdot 83833695555158789$,
$p_{n+3} = 1425172824437699414 = 2 \cdot 7 \cdot 11 \cdot 10603141 \cdot 872795051$,
$p_{n+4} = 1425172824437699415 = 3 \cdot 5 \cdot 13 \cdot 47 \cdot 155501672060851$
$p_{n+5} = 1425172824437699416 = 2^3 \cdot 29 \cdot 1400023 \cdot 4387775281$
…,
$p_{n+100} = 1425172824437699511 = 3^2 \cdot 1307 \cdot 121157257879597$,
…,
$p_{n+1000} = 1425172824437700411 = 3^5 \cdot 72883 \cdot 80470188419$,
…,
$p_{n+1473} = 1425172824437700884 = 2^2 \cdot 7 \cdot 13 \cdot 251 \cdot 2352079 \cdot 6631939$,
$p_{n+1474} = 1425172824437700885 = 3 \cdot 5 \cdot 599 \cdot 19657565689$,
$p_{n+1475} = 1425172824437700886 = 2 \cdot 712586412218850443$,
$p_{n+1476} = 1425172824437700887 = 1425172824437700887$,

seem to have random prime factorizations, not smooth prime factorizations as primordial integers do.

A few others reported examples for the prime gaps of very large probable primes also satisfy the estimate, in particular,

$$2254930 = p_{n+1} - p_n = 2\pi \log^2 p_n / \log\log p_n = 20587677614.2$$

for a 86853 decimal digits probable prime $p_n$, see [KR].

The statistical aspect of the prime gaps are derived from the Hardy-Littlewood conjecture for the number of primes $\pi_d(x) = \#\{ p \leq x : q = p + d \text{ is prime} \}$, see [KA] for recent work in this area. The numerical data for the frequency distribution of the prime gaps is given in [OK].

| $n$ | $d_n$ | $p_n$ | $U(n)$ | $n$ | $d_n$ | $p_n$ | $U(n)$ |
|---|---|---|---|---|---|---|---|
| 1 | 1 | 2 | -8.2 | 39 | 456 | 25056082087 | 1134.3 |
| 2 | 2 | 3 | 80.6 | 40 | 464 | 42652618343 | 1177.2 |
| 3 | 4 | 7 | 35.7 | 41 | 468 | 127976334671 | 1267.8 |
| 4 | 6 | 23 | 54.1 | 42 | 474 | 182226896239 | 1297.6 |
| 5 | 8 | 89 | 84.3 | 43 | 486 | 241160624143 | 1321.4 |
| 6 | 14 | 113 | 90.4 | 44 | 490 | 297501075799 | 1339.4 |
| 7 | 18 | 523 | 134.2 | 45 | 500 | 303371455241 | 1341.1 |
| 8 | 20 | 887 | 151.2 | 46 | 514 | 304599508537 | 1341.4 |
| 9 | 22 | 1129 | 159.2 | 47 | 516 | 416608695821 | 1368.5 |
| 10 | 34 | 1327 | 164.7 | 48 | 532 | 461690510011 | 1377.4 |
| 11 | 36 | 9551 | 238.2 | 49 | 534 | 614487453523 | 1402.4 |
| 12 | 44 | 15683 | 258.5 | 50 | 540 | 738832927927 | 1418.6 |
| 13 | 52 | 19609 | 267.9 | 51 | 582 | 1346294310749 | 1471.9 |
| 14 | 72 | 31397 | 288.2 | 52 | 588 | 1408695493609 | 1475.9 |
| 15 | 86 | 155921 | 362.0 | 53 | 602 | 1968188556461 | 1506.1 |
| 16 | 96 | 360653 | 403.6 | 54 | 652 | 2614941710599 | 1531.9 |
| 17 | 112 | 370261 | 404.9 | 55 | 674 | 7177162611713 | 1625.2 |
| 18 | 114 | 492113 | 419.5 | 56 | 716 | 13829048559701 | 1687.1 |
| 19 | 118 | 1349533 | 472.9 | 57 | 766 | 19581334192423 | 1720.3 |
| 20 | 132 | 1357201 | 473.2 | 58 | 778 | 42842283925351 | 1796.2 |
| 21 | 148 | 2010733 | 494.8 | 59 | 804 | 90874329411493 | 1870.4 |
| 22 | 154 | 4652353 | 542.2 | 60 | 806 | 171231342420521 | 1934.0 |
| 23 | 180 | 17051707 | 619.4 | 61 | 906 | 218209405436543 | 1958.6 |
| 24 | 210 | 20831323 | 631.7 | 62 | 916 | 1189459969825483 | 2134.4 |
| 25 | 220 | 47326693 | 683.3 | 63 | 924 | 1686994940955803 | 2171.4 |
| 26 | 222 | 122164747 | 745.0 | 64 | 1132 | 1693182318746371 | 2171.8 |
| 27 | 234 | 189695659 | 774.5 | 65 | 1184 | 4384154784554159 | 2274.2 |
| 28 | 248 | 191912783 | 775.2 | 66 | 1198 | 55350776431903243 | 2557.1 |
| 29 | 250 | 387096133 | 823.2 | 67 | 1220 | 80873624627234849 | 2600.7 |
| 30 | 282 | 436273009 | 831.5 | 68 | 1224 | 203986478517455989 | 2708.4 |
| 31 | 288 | 1294268491 | 908.8 | 69 | 1248 | 218034721194214273 | 2716.2 |
| 32 | 292 | 1453168141 | 917.2 | 70 | 1272 | 305405826521087869 | 2756.0 |
| 33 | 320 | 2300942549 | 950.8 | 71 | 1328 | 352521223451364323 | 2773.0 |
| 34 | 336 | 3842610773 | 989.0 | 72 | 1356 | 401429925999153707 | 2788.4 |
| 35 | 354 | 4302407359 | 997.6 | 73 | 1370 | 418032645936712127 | 2793.3 |
| 36 | 382 | 10726904659 | 1067.5 | 74 | 1442 | 804212830686677669 | 2871.7 |
| 37 | 384 | 20678048297 | 1119.1 | 75 | 1476 | 1425172824437699411 | 2941.0 |
| 38 | 394 | 22367084959 | 1125.3 | | | | |